\documentstyle{article}

\begin{document}
\begin{center}
{\Large\bf REVIEW OF MATHEMATICAL TECHNIQUES APPLICABLE IN ASTROPHYSICAL REACTION RATE THEORY}\\[1cm]
{\large A.M. MATHAI}\\ 
Department of Mathematics and Statistics, McGill University,\\
805 Sherbooke Street West, Montreal, CANADA H3A 2K6\\[0.5cm]
{\large H.J. HAUBOLD}\\ 
Office for Outer Space Affairs, United Nations,\\ 
P.O. Box 500, A-1400 Vienna, AUSTRIA\\[0.5cm]
\end{center}
\noindent
{\bf Abstract.} An overview is presented on statistical techniques for the analytic evaluation of integrals for non-resonant, non-resonant depleted, non-resonant cut-off, non-resonant sccreened, and resonant thermonuclear reaction rates. The techniques are based on statistical distribution theory and the theory of Meijer's G-function and Fox's H-function. The implementation of Meijer's G-function in Mathematica constituts an additional utility for analytic manipulations and numerical computation of thermonuclear reaction rate integrals. Recent results in the astrophysical literature related to the use of analytic thermonuclear reaction rates are incorporated.     
\section{Introduction}
Understanding the methods of evaluation of thermonuclear reaction rates is one of the most important goals of research in the field of stellar and cosmological nucleosynthesis. Practically all applications of fusion plasmas are controlled in some way or another by the theory of thermonuclear reaction rates under specific circumstances. After several decades of effort, a systematic and complete theory of thermonuclear reaction rates has been constructed (Lang, 1999). One of the basic ideas in this regard is that the motor of irreversibility and dissipation is the existence of reactions between individual nuclei. The latter produce a randomization of the energy and velocity distributions of particles. The effect of the reactions is balanced by the flow of the particles in a macroscopically inhomogeneous medium. As a result of this balance, the system reaches a quasi-stationary state close to equilibrium, in which steady fluxes of matter, energy, and momentum are present. The main ideas in the following are coming from statistical distribution theory and the theory of generalized special functions, mainly in the categories of Meijer's G-function and Fox's H-function of scalar, vector, and matrix arguments (Mathai, 1993; Aslam Chaudhry and Zubair, 2002). For an overview on the appliaction and historical background of such functions for particle distributions, see Hegyi (1999). An interesting application of them in the formation of structure in the universe, see Adler and Buchert (1999) and Buchert at al. (1999). 
A fusion of mathematical and statistical techniques enabled us to evaluate thermonuclear reaction rate integrals in explicit closed forms. Some of the techniques which are used will be summarized here. In order to explain the ideas we will start with the evaluation of an integral over a real scalar variable first. Let
\begin{equation}
I(z;p,n,m)=p\int^\infty_0e^{-pt}t^{-n\rho}e^{-zt^{-n/m}}dt
\end{equation}
for $\Re(p)>0, \Re(z)>0,$ n, m positive integers, where $\Re(.)$ denotes the real part of (.). A particular case of this for $n\rho =-\nu, p=1, n=1, m=2,$
\begin{equation}
I(z;1,1,2)=\int^\infty_0 e^{-t}t^\nu e^{-zt^{-1/2}}dt
\end{equation}
is a thermonuclear function associated with equilibrium distributions in reaction rate theory under Maxwell-Boltzmannian approach. As can be seen from (1)that the usual mathematical techniques fail to obtain a closed-form representation of the basic integral in (1).
\section{Statistical Techniques}
Certain special functions are related to particular probability laws governing products of independent exponential variables. Such laws can be related to the underlying physical processes. The integrand in (1) is a product of integrable positive functions and hence it can be made into a product of statistical densities by normalizing them. Consider two statistically independent real scalar random variables x and y with the density functions $f_1(x)\geq 0, f_2(y) \geq 0$ for $0 < x<\infty, 0<y<\infty$ and $f_1(x)=0, f_2(y)=0 $ elsewhere. Let $u=xy$, the product of these random variables. Then from the transformation of the variables, $u=xy, v=x$, the density $g(u)$ of $u$ is given by
\begin{equation}
g(u)=\int^\infty_0 \frac{1}{v}f_1(v)f_2(u/v)dv.
\end{equation}
The integral in (3) can be made equivalent to the integral in (1) by suitably selecting $f_1$ and $f_2$. Let
\begin{equation}
v=pt, u=pz^{m/n}, f_1(t)=t^{1-n\rho}e^{-t},\;\;\mbox {and}\;\; f_2(t)=e^{-t^{n/m}},
\end{equation}
excluding the normalizing constants. Then
$$
\int^\infty_0\frac{1}{v}f_1(v)f_2(u/v)dv=\int^\infty_0\frac{1}{t}e^{-pt}(pt)^{1-n\rho}e^{-zt^{-n/m}}dt.
$$
Hence
$$p^{p\rho}\int^\infty_0\frac{1}{t}e^{-pt}(pt)^{1-n\rho}e^{-zt^{n/m}}dt=p\int^\infty_0e^{pt}t{-n\rho}e^{-zt^{-n/m}}dt$$
which is exactly the integral to be evaluated in (1). We have identified the integral as the exact density of the product of two real scalar random variables, $u=xy$. Since this density is unique the idea is to evaluate this density through some other means. Notice that $u$ is a product of positive variables and hence the method of moments can be exploited profitably. 
Consider the $(s-1)$th moment of $u$, denoted by expected value of $u^{s-1}$. That is,
$$
E(u^{s-1})=E(x^{s-1})E(y^{s-1})
$$
due to statistical independence of $x$ and $y$. Let
$$
g_1(s)=E(x^{s-1})\;\;\mbox{and}\;\; g_2(s)=E(y^{s-1}).
$$
Then $g_1(s)$ and $g_2(s)$ are the Mellin transforms of $f_1$ and $f_2$ respectively. Then from the inverse Mellin transform the unique density of $u$ is available as
\begin{equation}
g(u)=\frac{1}{2\pi i}\int_Lu^{-s}g_1(s)g_2(s)ds, i=\sqrt{-1}.
\end{equation}
where $L$ is a suitable contour. But, excluding the normalizing constants,
\begin{eqnarray}
g_1(s) & = &\int^\infty_0t^{s-1}f_1(t)dt\nonumber\\
& = & \int^\infty_0t^{1-n\rho +s-1}e^{-t}dt\nonumber\\
 & = &\Gamma(1-n\rho+s)\;\; \mbox{for}\;\; \Re(1-n\rho+s)>0,\label{line3}.
\end{eqnarray} 
and
\begin{eqnarray}
g_2(s) & = &\int^\infty_0t^{s-1}f_2(t)dt\nonumber\\
& = & \int^\infty_0t^{s-1}e^{-t^{n/m}}dt\nonumber\\
& = & \frac{m}{n}\Gamma(ms/n)\;\; \mbox{for}\;\;\Re(s)>0\label{line3}.
\end{eqnarray}
Then
\begin{equation}
\frac{1}{2\pi i}\int_Lu^{-s}g_1(s)g_2(s)ds=\frac{1}{2\pi i}\int_L(m/n)\Gamma(ms/n)\Gamma(1-n\rho+s)(pz^{m/n})^{-s}ds.
\end{equation}
Therefore
\begin{eqnarray}
p & \int^\infty_0 &  e  ^{-pt}t^{-n\rho}e^{-zt^{-n/m}}dt\nonumber\\
 & = &  p^{n\rho}\int^\infty_0\frac{1}{t}e^{-pt}(pt)^{1-n\rho} e^{-zt^{-n/m}}dt\nonumber\\
 & = & p^{n\rho}(m/n)\frac{1}{2\pi i}\int_L\Gamma(ms/n)\Gamma(1-n\rho+s)(pz^{m/n})^{-s}ds\nonumber\\
 & = & mp^{n\rho}\frac{1}{2\pi i}\int_{L_1}\Gamma(ms)\Gamma(1-n\rho +ns)(z^mp^n)^{-s}ds
\end{eqnarray}
by replacing $s/n$ by $s$. Our aim now is to evaluate the contour integral on the right side of (9) explicitly into computable forms. The contour integral in (9) can be written as an H-function which can then be reduced to a G-function since $m$ and $n$ are positive integers.
\section{G- and H-Functions}
For the sake of completeness we will define G- and H-functions here. For the theory and applications of G-functions see Mathai (1993) and their relation to incomplete gamma functions (Aslam Chaudhry and Zubair, 2002). Let $\phi(s)$ and $\psi(s)$ be the following gamma products.
$$\psi(s)=\frac{\left\{\Pi^m_{j=1}\Gamma(b_j+\beta_js)\right\}\left\{\Pi^n_{j=1}\Gamma(1-a_j-a_js)\right\}}{\left\{\Pi^q_{j=m+1}\Gamma(1-b_j-\beta_js)\right\}\left\{\Pi^p_{j=n+1}\Gamma(a_j+\alpha_js)\right\}}$$
and
$$\phi(s)=\psi(s)\;\;\mbox{with}\;\; \alpha_j=1, j=1,\ldots,p, \beta_k=1,k=1,\ldots, q.$$ 
Then Fox's H-function, denoted by
$$H^{m,n}_{p,q}(z)=H^{m,n}_{p,q}[z|^{(a_1, \alpha_1),\ldots,(a_p, \alpha_p)}_{(b_1, \beta_1),\ldots,(b_q, \beta_q)}]
$$
is defined as 
\begin{equation}
H_{p,q}^{m,n}(z)=\int_{L_2}\psi(s)z^{-s}ds
\end{equation}
where $L_2$ is a suitable contour separating the poles of the gammas $\Gamma(1-a_j-\alpha_j s), j=1,\ldots , n$ from those of $\Gamma(b_j+\beta_js), j=1,\ldots,m.$ 
The theory and applications of H-functions are available from Mathai and Saxena (1978). In (10) when $\alpha_j=1, j=1,\ldots, p\;\; \mbox{and}\;\; \beta_k=1, k=1, \ldots, q$ we have Meijer's G-function denoted by
\begin{eqnarray}
G^{m,n}_{p,q} & = & G^{m,n}_{p,q}[z|^{a_1,\ldots,a_p}_{b_1,\ldots,b_q}]\nonumber\\
& = & \int_{L_3}\phi(s)z^{-s}ds
\end{eqnarray}
where $L_3$ is a suitable contour. Various types of contours, existence conditions and properties of the G-function are available from Mathai (1993). \par
Now, comparing (9) and (10) our starting integral is evaluated as follows
\begin{eqnarray}
I(z;p,n,m)& = & mp^{n\rho}\frac{1}{2\pi i}\int_{L_1}\Gamma(ms)\Gamma(1-n\rho+ns)(z^mp^n)^{-s}ds\nonumber\\
& = & mp^{n\rho}H^{2,0}_{0,2}[(z^mp^n)|_{(0,m),(1-n\rho,n)}].
\end{eqnarray}
The H-function in (12) can be reduced to a G-function which can again be reduced to computable series forms. For this purpose we expand the gammas in the integrand in (12) by using the multiplication formula for gamma functions, namely,
\begin{equation}
\Gamma(mz)=(2\pi)^{\frac{(1-m)}{2}}m^{mz-\frac{1}{2}}\Gamma(z)\Gamma(z+\frac{1}{m})\ldots\Gamma(z+\frac{m-1}{m})
\end{equation}
$$ m  = 1,2,\ldots $$
Expanding $\Gamma(ms)\Gamma(1-n\rho+ns)$ by using (13) we have
\begin{eqnarray*}
mp^{n\rho}\Gamma(ms)\Gamma(1-n\rho +ns)& = & (2\pi)^{\frac{1}{2}(2-m-n)}p^{n\rho}m^{\frac{1}{2}}n^{-n\rho+\frac{1}{2}}(m^mn^n)^s\nonumber\\
& \times & \Gamma(s)\Gamma(s+\frac{1}{m})\ldots\Gamma(s+\frac{m-1}{m})\nonumber\\
& \times & \Gamma(s+\frac{1-n\rho}{n})\ldots \Gamma(s+\frac{n-n\rho}{n})
\end{eqnarray*}
Substituting these back and writing as a G-function we have
\begin{eqnarray}
I(z;p,n,m) &=&p^{n\rho}(2\pi)^{\frac{1}{2}(2-m-n)}m^{\frac{1}{2}}n^{-n\rho+\frac{1}{2}}\nonumber\\
& \times & G^{m+n,0}_{0,m+n}[\frac{z^mp^n}{m^mn^n}|_{0,\frac{1}{m},\ldots,\frac{m-1}{m}, \frac{1-n\rho}{n},\ldots,\frac{n-n\rho}{n}}].
\end{eqnarray}
\section{Non-Resonant Thermonuclear Reaction Rate}
The Maxwell-Boltzmannian form of the collision probability integral for non-resonant thermonuclear reactions is (Lang 1999; Haubold and Mathai, 1984; Hussein and Pato 1997; Bergstroem et al. 1999; Ueda et al. 2000) 
\begin{eqnarray}
I_1& = & I(z;1,1,2)=\int^\infty_0 y^\nu e^{-y}e^{-zy^{-1/2}}dy\nonumber\\
& = & \pi^{-\frac{1}{2}}G^{3,0}_{0,3}[\frac{z^2}{4}|_{0,\frac{1}{2}, 1+\nu}].
\end{eqnarray}

In stellar fusion plasmas the energies of the moving nuclei are assumed to be described by a Maxwell-Boltzmann distribution, $~E exp(-E/kt)$, where $T$ is the local temperature and $k$ the Boltzmann constant. Folding the cross section of a nuclear reaction, $\sigma(E)$, with this energy (or velocity) distribution leads to the nuclear reaction rate per pair of nuclei:
\[<\sigma v> = (8/\pi \mu)^{1/2}(kT)^{-3/2} \int_0^\infty \sigma(E)exp(-E/kT)dE,\]
where $v$ is the relative velocity of the pair of nuclei, $E$ is the center-of-mass energy, and $\mu =m_1m_2/(m_1+m_2)$ is the reduced mass of the entrance channel of the reaction. In order to cover the different evolution phases of the stars, i.e. from main sequence stars to supernovae, one must know the reaction rates over a wide range of temperatures, which in turn requires the availability of $\sigma(E)$ data over a wide range of energies. I is the challenge to the experimentalist to make precise $\sigma(E)$ measurements over a wide range of energies. For the class of charged-particle-induced reactions, there is a repulsive Coulomb barrier in the entrance channel of height $E_c=Z_1Z_2e^2/r$, where $Z_1$ and $Z_2$ are the integral nuclear charges of the interacting particles, e is the unit of electric charge, and r is the nuclear interaction radius. Due to the tunneling effect through the Coulomb barrier, $\sigma(E)$ drops nearly exponentially with decreasing energy:
\[\sigma(E)=S(E)E^{-1}exp(-2\pi\eta),\]
where $\eta=Z_1Z_2e^2/hv$ is the Sommerfeld parameter, $h$ is the Planck constant. The function $S(E)$ contains all the strictly nuclear effects, and is referred to as the astrophysical $S(E)$ factor. If the above equation for $\sigma(E)$ is inserted in the above equation for the nuclear reaction rate $<\sigma v>$, one obtains
\[<\sigma v>=(8/\pi\mu)^{1/2}(kT)^{-3/2} \int_0^\infty S(E)exp(-E/kT-b/E^{1/2})dE,\]
with $b=2(2\mu)^{1/2}\pi^2e^2Z_1Z_2/h.$ Since for nonresonant reactions $S(E)$ varies slowly with energy, the steep energy dependence of the integrand in the equation for $<\sigma v>$ is governed by the exponential term, which is characterised by the peak near an energy E0 that is usually much larger than $kT$, the mean termal energy of the fusion plasma. The peak is referred to as the Gamow peak. For a constant $S(E)$ value over the energy region of the peak, one finds $E_0=(bkT/2)^{2/3}.$ This is the effective mean energy for a givern reaction at a given temperature. If one approximates the peak by a Gaussian function, one finds an effective width $\delta=4(E_0kT)^{1/2}/3^{1/2}.$ In the following, this approximation is not made and the respective integrals are analytically represented, beginning with  

\begin{eqnarray*}
<\sigma v> & = &(\frac{8}{\pi\mu})^{\frac{1}{2}}\sum^2_{\nu=0}\frac{1}{(kT)^{-\mu+\frac{1}{2}}}\frac{S^{(\mu)}(0)}{\mu!}\nonumber\\
& \times & \int^\infty_0 e^{-y}y^\nu e^{-zy^{-1/2}}dy
\end{eqnarray*}
where $S^(\mu)$ denotes the $\mu$-the derivative. 
The G-function in (15) can be expressed as a computable power series as well as in closed-forms by using residue calculus. Writing the G-function in (15) as a Mellin-Barnes integral we have
\begin{equation}
G^{3,0}_{0,3}[\frac{z^2}{4}|_{0,\frac{1}{2},1+\nu}]=\frac{1}{2\pi i}\int\Gamma(s)\Gamma(s+1/2)\Gamma(1+\nu+s)(z^2/4)^{-s}ds.
\end{equation}
{\bf Case (1):} $\nu\neq \pm\frac{\lambda}{2}, \lambda = 0,1,2,\ldots$
In this case the poles of the integrand are simple and the poles are at the points
$$s=0,-1,-2,\ldots; s=-\frac{1}{2}, -\frac{1}{2}-1,\ldots; s= -1-\nu,-2-\nu,\ldots$$
\begin{eqnarray}
\sum^\infty_{r=0}& \frac{(-1)^r}{r!} & \Gamma(\frac{1}{2}-r)\Gamma(1+\nu-r)(\frac{z^2}{4})^r\nonumber\\
& = & \Gamma(\frac{1}{2})\Gamma(1+\nu)_0F_2(;\frac{1}{2},-\nu;-\frac{z^2}{4})
\end{eqnarray}
where, in general, $_pF_q$ denotes a general hypergeometric function. The sum of the residues at $s=-\frac{1}{2}, -\frac{1}{2} -1, \ldots$ is
\begin{eqnarray}
\sum^\infty_{r=0} & \frac{(-1)^r}{r!}& \Gamma(-\frac{1}{2}-r)\Gamma(\frac{1}{2}+\nu-r)(\frac{z^2}{4})^{\frac{1}{2}+r}\nonumber\\
&=& \Gamma(-\frac{1}{2})\Gamma(\frac{1}{2}+\nu)(\frac{z^2}{4})^{\frac{1}{2}}\;\;_0F_2(;\frac{3}{2},\frac{1}{2}-\nu;-\frac{z^2}{4}). 
\end{eqnarray}
The sum of the residues at $s=-1-\nu, -1-\nu-1,\ldots$ is
\begin{eqnarray}
\sum^\infty_{r=0} & \frac{(-1)^r}{r!}& \Gamma(-\frac{1}{2}-\nu -r)\Gamma(-1-\nu -r)(\frac{z^2}{4})^{1+\nu +r}\nonumber\\
& = & \Gamma(-1-\nu)\Gamma(-\frac{1}{2}-\nu)(\frac{z^2}{4})^{1+\nu}\;\;_0 F_2(;\nu+2,\nu+\frac{3}{2};-\frac{z^2}{4})
\end{eqnarray}
Then from (15) to (19) we have
\begin{eqnarray}
I(z;1,1,2) & = & \int_0^\infty y^\nu e^{-y} e^{-zy^{1/2}} dy\nonumber\\
& = & \pi^{-\frac{1}{2}}G^{3,0}_{0,3}[\frac{z^2}{4}|_{0,\frac{1}{2}, 1+\nu}]\nonumber\\
& = & \Gamma(1+\nu) _0 F_2(;\frac{1}{2},-\nu; -\frac{z^2}{4})\nonumber\\
& -& 2\Gamma(\frac{1}{2}+\nu)(\frac{z^2}{4})^{\frac{1}{2}}\;\;_0 F_2(; \frac{3}{2}; \frac{1}{2}-\nu; -\frac{z^2}{4})\nonumber\\
& + & \frac{\Gamma(-1-\nu)\Gamma(-\frac{1}{2}-\nu)}{\gamma(\frac{1}{2})}(\frac{z^2}{4})^{1+\nu}\nonumber\\
& \times & _0 F_2(; \nu+2,\nu+\frac{3}{2}; -\frac{z^2}{4})
\end{eqnarray}
for $\nu\neq \pm\frac{\lambda}{2}, \lambda = 0,1,2,\ldots$ When $\nu$ is a positive integer the poles at $s=-1-\\nu, -1-\nu-1,\ldots$ are of order 2 each. Hence the corresponding sum of residues can be written in terms of psi functions. Similarly when $\nu$ is a negative integer, positive or negative half integer the corresponding  sums will contain psi functions. Details of the computations can be seen from Haubold and Mathai (1984) and their application by Hussein and Pato (1997), Bergstroem et al. (1999), and Ueda et al. (2000).
\section{Modified Non-Resonant Thermonuclear Reaction Rate: Depletion}
With deviations from the Maxwell-Boltzmann velocity distribution of nuclei in the fusion plasma, a modification which results in the depletion of the tail is introduced in Haubold and Mathai (1986a); see also Kaniadakis et al. (1997), Kaniadakis et al. (1998), and Coraddu et al. (1999). In this case the collision probability integral will be of the following form:
\begin{equation}
I_2=\int_0^\infty y^\nu e^{-y}e^{-zy^{-1/2}}dy.
\end{equation}
We consider a general integral of this type. Let
\begin{equation}
I(z;\delta,a,b,m,n)=\int_0^\infty t^\rho e^{-at-bt^\delta-zt^{-n/m}}dt.
\end{equation}
Expanding $e^{-bt^\delta}$ and then with the help of (14) one can represent (22) in terms of a G-function as follows:
\begin{eqnarray}
I(z;\delta,a,b,m,n) & = &\sum^\infty_{k=0}\frac{(-b)^k}{k!}a^{-(\rho+k\delta+1)}(2\pi)^{\frac{1}{2}(2-m-n)}m^{\frac{1}{2}}n^{\frac{1}{2}+\rho+k\delta}\nonumber\\
& \times & G^{m+n,0}_{0,m+n}[\frac{z^ma^n}{m^mn^n}|_{0,\frac{1}{m},\ldots,\frac{m-1}{m}, \frac{1+\rho+k\delta}{n},\ldots,\frac{n+\rho+k\delta}{n}}]
\end{eqnarray}
for $\Re(z)>0,\Re(a)>0,\Re(b)>0,m,n=1,2,\ldots.$ the case in (21) is for $a=1,b = 1,n = 1,m=2$. With $\nu=\rho+k\delta$ the G-function in (21) corresponds to that in (15). When $\delta$ is irrational and $\rho$ is rational, the poles of the integrand will be simple and the G-function is available in terms of hypergeometric functions. Other situations will involve psi functions. From the asymptotic behavior of the G-function, see for example Mathai (1993), one can write the integral in (21), for large values of $z$ as follows:
\begin{equation}
I_2\approx \pi^{\frac{1}{2}}(\beta/3)^{\frac{2\nu+1}{2}}e^{-\beta-(\beta/3)^\delta}, \beta=3(z/2)^{2/3}.
\end{equation}
\section{Modified Non-Resonant Thermonuclear Reaction Rate: Cut-Off}
Another modification can be made by acut-off of the high-energy tail of the Maxwell-Boltzmann distribution. In this case the collision probability integral to be evaluated is of the form
\begin{equation}
I_3=\int^d_0t^{-\rho}e^{-at-zt^{-1/2}}dt, d<\infty,
\end{equation}
see the details from Haubold and Mathai (1986b). We will consider a general integral of the form
\begin{equation}
I(z;d,a,\rho,n,m)=\int^d_0 y^{-n\rho} e^{-ay-zy^{-n/m}}dy, d<\infty.
\end{equation}
In order to evaluate (26) explicitly we will use statistical techniques as discussed earlier. Let $x$ and $y$ be two statistically independent real random variables having the densities $c_1f_1(x), 0<x<d$ and $c_2f_2(y), 0<y<\infty$ with $f_1(x)$ and $f_2(y)$ equal to zero elsewhere, where $c_1$ and $c_2$ are normalizing constants. Then taking
$$f_1(x)=x^{-n\rho+1}e^{-ax}\;\; \mbox{and}\;\; f_2(y)=e^{-y^{n/m}}$$
and proceeding as before one has the following result:
\begin{eqnarray}
I(z;d,a,\rho,n,m) & =& \int^d_0t^{-n\rho}e^{-at-zt^{-n/m}} dt\nonumber\\
& = & m^{\frac{1}{2}}n^{-1}(2\pi)^{(1-m)/2}d^{-n\rho+1}\sum^\infty_{r=0}\frac{(-ad)^r}{r!}\nonumber\\
& \times & G^{m+n,0}_{n,m+n}\left[\frac{z^m}{d^n m^m}|^{-\rho+\frac{r+2+j-1}{n},j=1,\ldots,n}_{-\rho+\frac{r+1+j-1}{n}, j=1,\ldots,n;\frac{j-1}{m},j=1,\ldots,m}\right]
\end{eqnarray}
for $\Re(z)>0,>0,\Re(a)>0.$ Then
\begin{eqnarray}
I_3 & = & \int^d_0t^{-\rho} e^{-at-zt^{-/2}}dt\nonumber\\
& = & d^{-\rho+1}\pi^{-\frac{1}{2}}\sum^\infty_{r=0}\frac{(-ad)^r}{r!}\nonumber\\
& \times & G^{3,0}_{1,3}\left[\frac{z^2}{4d}|^{-\rho+r+2}_{-\rho+r+1,0,\frac{1}{2}}\right].
\end{eqnarray}
For large values of $z$ the G-function behaves like $\pi^{\frac{1}{2}}x^{-\frac{1}{2}}e^{-2x^{\frac{1}{2}}}, x=\frac{z^2}{4d}$, see for example Mathai (1993). Then for large values of $z$,
\begin{equation}
I_3\approx d^{-\rho+1}\left(\frac{z^2}{4d}\right)^{-\frac{1}{2}}\;\;e^{-ad-2(z^2/4d)^\frac{1}{2}}.
\end{equation}
Explicit series forms can be obtained for various values of the parameters with the help of residue calculus. For example, for $\Re(z)>0,d>0,\Re(a)>0,-\rho+r+1\neq\mu, \mu=0,1,\ldots$
\begin{eqnarray}
I_3 &=&\int_0^d t^{-\rho}e^{-at-zt^{-1/2}}dt\nonumber\\
&=&\pi^{-\frac{1}{2}}d^{-\rho+1}\sum^\infty_{r=0}\frac{(-ad)^r}{r!}\nonumber\\
&\times &\left\{ \sum^\infty_{\nu=0,\nu\neq\nu}\frac{(-1)^\nu\Gamma(\frac{1}{2}-\nu)}{\nu!(-\rho+r+1-\nu)}(\frac{z^2}{4d})^\nu\right.\nonumber\\
& + & \sum^\infty_{\nu=0}\frac{(-1)^\nu \Gamma(-\frac{1}{2}-\nu)}{\nu!(-\rho +r-\nu+\frac{1}{2})}(\frac{z^2}{4d})^{\nu+
\frac{1}{2}}\nonumber\\
&+&\left.(\frac{z^2}{4d})^\mu\left[-ln(\frac{z^2}{4d})+A\right]B\right\},
\end{eqnarray}
where
$$A=\psi(\mu+1)+\psi(-\mu+\frac{1}{2}), \psi(z)=\frac{d}{dz}ln\Gamma(z),$$
$$B=\frac{(-1)^\mu}{\mu!}\Gamma (-\mu+\frac{1}{2}).$$
Other cases and further details may be seen from Haubold and Mathai (1986b).
\section{Screened Non-Resonant Thermonuclear Reaction Rate}
Taking into account plasma corrections to fusion processes due to a static or dynamic potential, see the details from Haubold and Mathai (1986c), Shaviv and Shaviv (2001), Wierling et al. (2001), Bahcall et al. (2002), the collision probability integral to be evaluated is of the form
\begin{equation}
I_4=\int^\infty_0y^n e^{-y-z(y+t)^{-1/2}}dy.
\end{equation}We will consider a more general integral of the following form:
\begin{equation}
I(z;t;a,\nu,n,m)=\int^\infty_0 y^\nu e^{-ay-z(y+t)^{-n/m}}dy.
\end{equation}
This can be evaluated with the help of (14) and (27). Let
\begin{eqnarray*}
N_1(z)& = & N_1(z;a,\nu, n,m)\\
& = & \int^\infty_0 y^\nu e^{-ay-zy^{-n/m}}dy.
\end{eqnarray*}
Then from (14)
\begin{eqnarray}
N_1(z) = a^{-(\nu+1)}(2\pi)^{\frac{1}{2}(2-m-n)}m^{\frac{1}{2}}n^{\nu+\frac{1}{2}}\nonumber\\
\times  G^{m+n,0}_{0,m+n}\left[\frac{z^ma^n}{m^mn^n}|_{0,\frac{1}{m},\ldots,\frac{m-1}{m},\frac{n+\nu}{n},\ldots,\frac{n+\nu}{n}}\right]
\end{eqnarray}
for $\Re(a)>0,\Re(z)>0, m,n$ positive integers. Let
\begin{eqnarray*}
N_2(z)&= & N_2(z;d,a,\nu,n,m)\\
& = & \int_0^dy^\nu e^{-ay-zy^{-n/m}}dy.
\end{eqnarray*}
Then from (27)
\begin{eqnarray}
N_2(z)&=&(2\pi)^{(1-m)/2}m^{\frac{1}{2}}n^{-1}d^{\nu+1}\sum^\infty_{r=0}\frac{(-ad)^r}{r!}\nonumber\\
&\times&G^{m+n,0}_{n,m+n}[\frac{z^m}{d^nm^m}|^{\frac{\nu+r+2+j-1}{n},j=1,\ldots,n}_{\frac{\nu+r+1+j-1}{n},j=1,\ldots,n,\frac{j-1}{m},j=1,\ldots,m}]
\end{eqnarray}
for $\Re(z)>0,d>0,\Re(a)>0,m,n$ positive integers. Then by a change of variables and rewriting the integral we have
\begin{eqnarray}
I_4 &=&\int_0^\infty y^\nu e^{-ay-z(y+t)^{-n/m}}dy \nonumber\\
&=& t^{\nu+1}e^{a_1}\sum^\nu_{r=0}(^\nu_r)(-1)^r\nonumber\\
&\times&[N_1(z_1;a_1,\nu-r,n,m)-N_2(z_1;1,a_1,\nu-r,n,m)]
\end{eqnarray}
where
$$a_1=at, z_1=zt,(^\nu_r)=\frac{\nu!}{r!(\nu-r)!}, 0!=1,$$
and $N_1(.)$ and $N_2(.)$ are defined in (33) and (34), respectively.

\section{Resonant Thermonuclear Reaction Rate}
When a resonance occurs in the low-energy range of nuclear reactions the collision probability integral can have the following form, for details see Haubold and Mathai (1986d).
\begin{equation}
R(q,a,b,g)=\int^\infty_0 t^\nu\frac{e^{-ay-qy^{-1/2}}}{(b-y)^2+g^2}dy.
\end{equation}
We will consider a more general integral of the following type:
\begin{equation}
R_1(q,a,b,g,\nu,n,m)=\int^\infty_0 t^\nu\frac{e^{-at-qt^{-n/m}}}{(b-t)^2+g^2}dt.
\end{equation}
The method to be employed here for evaluating this integral is to replace the denominator by an equivalent integral. That is, for $g^2>0$,
\begin{equation}
\frac{1}{(b-t)^2+g^2}=\int^\infty_0 e^{-[(b-t)^2+g^2]x}dx.
\end{equation}
Then expand
$$e^{-x(b-t)^2}=\sum^\infty_{k=0}\frac{(-1)^k}{k!}x^k\sum^{2k}_{k_1=0}(^{2k}_{k_1}(-1)^{k_1} b^{2k-k_1}t^{k_1}.$$
Substituting back in (37), integrating out $t$ first and then integrating over $x$ one has the following result:
\begin{eqnarray}
\int^\infty_0 &t^{\nu+k_1}& e^{-at-qt^{-n/m}}dt\nonumber\\
&=&a^{-(\nu+1+k_1)}(2\pi)^{\frac{1}{2}(2-n-m)}m^{\frac{1}{2}}n^{\nu+k_1+\frac{1}{2}}\nonumber\\
&\times& G^{m+n,0}_{0,m+n}\left[\frac{q^m a^n}{m^m n^n}|_{0,\frac{1}{m},\ldots,\frac{m-1}{m}, \frac{1+\nu+k_1}{n},\ldots,\frac{n+\nu+k_1}{n}}\right]
\end{eqnarray}
and
\begin{equation}
\int^\infty_0x^ke^{-g^2x}dx=\frac{k!}{(g^2)^{k+1}}.
\end{equation}
From (39) and (40)
\begin{eqnarray}
R_1&(q,&a,b,g,\nu,n,m)\nonumber\\
&=&\sum^\infty_{k=0}\frac{(-1)^k}{(g^2)^{k+1}}\sum^{2k}_{k_1=0}\left(^{2k}_{k_1}\right)(-1)^{k_1} b^{2k-k_1}\nonumber\\
&\times& a^{-(\nu+1+k_1)}(2\pi)^{\frac{1}{2}(2-n-m)}m^{\frac{1}{2}}n^{\nu+k_1+\frac{1}{2}}\nonumber\\
&\times& G^{m+n,0}_{0,m+n}\left[\frac{q^m a^n}{m^m n^n}|_{0,\frac{1}{m},\ldots,\frac{m-1}{m},\frac{1+\nu+{k_1}}{n},\ldots,\frac{n+\nu+{k_1}}{n}}\right].
\end{eqnarray}
Then for $n=1, m=2, \nu=0,$
\begin{eqnarray*}
R(q,a,b,g)&=&\frac{1}{\pi^{\frac{1}{2}}g^2a}\sum^\infty_{k=0}\frac{(-1)^k}{(g^2)^k}\nonumber\\
&\times&\sum^{2k}_{k_1=0}\left(^{2k}_{k_1}\right)\frac{(-1)^{k_1}}{a^{k_1}}b^{2k-k_1}\nonumber\\
&\times& G^{3,0}_{0,3}\left[\frac{q^2a}{4}|_{0,\frac{1}{2},k_1}\right]
\end{eqnarray*}
for $\frac{(b-\frac{\nu}{a})^2}{g^2}<1, \nu=(\frac{q^2a}{4})^{\frac{1}{3}}$.
\section{Computations}
For computational purposes we will consider the four basic integrals associated with the cases: non-resonant reactions, non-resonant "cut-off" reactions, non-resonant screened reactions, and non-resonant "depleted" reactions. Let
\begin{eqnarray}
J_1(z,\nu)&=&\int^\infty_0 y^\nu e^{-y-zy^{-1/2}}dy\nonumber\\
J_2(z,d,\nu)&=&\int^d_0 y^\nu e^{-y-zy^{-1/2}}dy\nonumber\\
J_3(z,t,\nu)&=&\int^\infty_0 y^\nu e^{-y-z(y+t)^{-\frac{1}{2}}}dy\nonumber\\
J_4(z,\delta,b,\nu)&=& \int^\infty_0 y^\nu e^{-y-by^\delta-zy^{-1/2}}dy.
\end{eqnarray}
The exact expressions for these are given in (15), (28), (35), (21) respectively. The symbolic evaluation of all these integrals can not yet be achieved with Mathematica (Wolfram, 1999). Those inetrals that involve no singularities are done by taking limits of the indefinite integrals. The definite versions of the integrals are done using the Marichev-Adamchik Mellin transform methods (Adamchik 1996). The integration results are initially expressed in terms of Meijer's G-function, which are subsequently converted into hypergeometric functions using Slater's theorem. The notation for Meijer's G-function, belonging to the implemented special functions of Mathematica, is
\begin{equation}
\mbox{MeijerG}[\left\{\left\{a_1,...,a_n\right\},\left\{a_{n+1},...,a_p\right\}\right\},\left\{\left\{b_1,...,b_m\right\},\left\{b_{m+1},...,b_q\right\}\right\},z].
\end{equation}
Analytic expressions for the following Meijer's G-functions are available on Wolfram Research's Mathematical Functions web page (Wolfram, 2002):
\begin{equation}
G\left\{m,n,p,q\right\}=G\left\{3,0,0,3\right\}=\mbox{http://functions.wolfram.com/07.34.03.0948.01,}
\end{equation} 
and
\begin{equation}
G\left\{m,n,p,q\right\}=G\left\{3,0,1,3\right\}=\mbox{http://functions.wolfram.com/07.34.03.0955.01.}
\end{equation}

Approximations for large values of $z$ can be worked out with the help of the asymptotic behavior of G-functions, see for example Mathai (1993). These are the following for $z$ very large:
\begin{eqnarray}
J_1& \approx & 2\left(\frac{\pi}{3}\right)^{\frac{1}{2}}\left(\frac{z^2}{4}\right)^{\frac{2\nu+1)}{6}} e^{-3(z^2/4)^\frac{1}{3}}\nonumber\\
J_2&\approx & d^{\nu+1}\left(\frac{z^2}{4d}\right)^{-\frac{1}{2}} e^{-d-2(z^2/4d)^{\frac{1}{2}}}\nonumber\\
J_3&\approx& 2\left(\frac{\pi}{3}\right)^{\frac{1}{2}}\left(\frac{z^2}{4}\right)^{\frac{1}{6}}\left[\left(\frac{z^2}{4}\right)^{\frac{1}{3}}-t\right]^\nu e^{t-3(z^2 /4)^{\frac{1}{3}}}\nonumber\\
J_4&\approx& 2\left(\frac{\pi}{3}\right)^{\frac{1}{2}}\left(\frac{z^2}{4}\right)^{\frac{2\nu+1}{6}}e^{-3(z^2/4)^{\frac{1}{3}}-b(z^2/4)^{\delta/3}}.
\end{eqnarray}
A number of exact and approximate graphs of these four integrals, for various values of the parameters, using the integration routines in Mathematica (Wolfram, 1999), are available from Anderson, Haubold and Mathai (1994). 
\section{A Generalization}
A mathematically interesting integral corresponding to (1) can be evaluated. Consider the integral
\begin{equation}
I=\int^\infty_0 e^{-pt}t^{\rho-1}e^{-zt^{-\gamma}}dt.
\end{equation}
Then take $f_1(x)=c_1 x^\rho e^{-px}, x>0, f_2(y)=c_2 e^{-y^\gamma}, \gamma>0, y>0$ and $f_1(x)=0,f_2(y)=0$ elsewhere, where $c_1$ and $c_2$ are normalizing constants. Then $u=xy=z^{1/\gamma}$ and from (3) one has the integral in (44) evaluated
as the following:
$$I=(\gamma p^\rho)^{-1}H^{2,0}_{0,2}[p z^{1/\gamma}|_{(p,1),(0,1/\gamma)}], 0<z<\infty,$$
where $H(.)$ is the H-function defined in (10). When $\gamma$ is rational the H-function can be rewritten in terms of a Meijer's G-function and then (44) can be evaluated in terms of the result given in (14). For specified values of $\gamma$ and $\rho$ one can obtain computable representations for the H-function. 
\bigskip
\begin{center}
{\bf References}
\end{center}
Adamchik, V.: 1996, Definite integration in Mathematica V3.0, {\it Mathematica\par
in Education and Research} {\bf 5}, no.3, 16-22.\par
\medskip
\noindent
Adler, S. and Buchert, T.: 1999, Lagrangian theory of structure formation\par
in pressure-supported cosmological fluids, {\it Astronomy and Astrophysics}\par
{\bf 343}, 317-324.\par
\medskip
\noindent
Anderson, W.J., Haubold, H.J., and Mathai, A.M.: 1994, Astrophysical\par
thermonuclear functions, {\it Astrophysics and Space Science} {\bf 214}, 49-70.\par
\medskip
\noindent
Aslam Chaudhry, M. and Zubair, S.M.: 2002, {\it On a Class of Incomplete Gamma\par
Functions with Applications}, Chapman and Hall/CRC, New York.\par
\medskip
\noindent
Bahcall, J.N., Brown, L.S., Gruzinov, A., and Sawyer, R.F.: 2002, The Salpeter\par
plasma correction for solar fusion reactions, {\it Astronomy and Astrophysics}\par
{\bf 383}, 291-295.\par
\medskip
\noindent
Bergstroem, L., Iguri, S., and Rubinstein, H.: 1999, Constraints on the variation\par
of the fine structure constant from big bang nucleosynthesis, {\it Physical Review}\par
{\bf D60}, 045005-1.\par
\medskip
\noindent
Buchert, T., Dominguez, A., and Perez-Mercader, J.: 1999, Extending the scope of\par
models for large-scale structure formation in the universe, {\it Astronomy and\par
Astrophysics} {\bf 349}, 343-353.\par
\medskip
\noindent
Coraddu, M., Kaniadakis, G., Lavagno, A., Lissia, M., Mezzorani, G., and\par
Quarati, P.: 1999, Thermal distributions in stellar plasmas, nuclear reactions\par
and solar neutrinos, {\it Brazilian Journal of Physics} {\bf 29}, 153-168.\par
\medskip
\noindent
Haubold, H.J. and Mathai, A.M.: 1984, On nuclear reaction rate theory,\par
{\it Annalen der Physik (Leipzig),}{\bf 41}(6), 380-396.\par
\medskip
\noindent
Haubold, H.J. and Mathai, A.M.: 1986a, Analytic representations of modified\par
non-resonant thermonuclear reaction rates, {\it Journal of Applied Mathematics\par
and Physics (ZAMP)} {\bf 37(5)}, 685-695.\par
\medskip
\noindent
Haubold, H.J. and Mathai, A.M.: 1986b, Analytic representations of\par
thermonuclear reaction rates, {\it Studies in Applied Mathematics} {\bf 75},\par
123-138.\par
\medskip
\noindent
Haubold, H.J. and Mathai, A.M.: 1986c, Analytic results for screened\par
non-resonant nuclear reaction rates, {\it Astrophysics and Space Science}\par
{\bf 127}, 45-53.\par
\medskip
\noindent
Haubold, H.J. and Mathai, A.M.: 1986d, The resonant thermonuclear reaction\par
rates, {\it Journal of Mathematical Physics} {\bf 27(8)}, 2203-2207.\par
\medskip
\noindent
Hegyi, S.: 1999, A powerful generalization of the NBD suggested by Peter\par
Carruthers, in {\it Correlations and Fluctuations '98 - From QCD to Particle\par
Interferometry: Proc. VIII Int. Workshop on Multiparticle Production},\par
Matrahaza, Hungary, 14-21 June 1998, Eds. T. Csoergo, S. Hegyi, G. Jancso,\par
and R.C. Hwa, World Scientific, pp. 272-286.\par
\medskip
\noindent
Hussein, M.S. and Pato, M.P.: 1997, Uniform expansion of the thermonuclear\par
reaction rate formula, {\it Brazilian Journal of Physics} {\bf 27}, no.3,\par
364-372.\par
\medskip
\noindent 
Kaniadakis, G., Lavagno, A., and Quarati, P.: 1997, Non-extensive statistics\par
and solar neutrinos, astro-ph/9701118.\par
\medskip
\noindent
Kaniadakis, G., Lavagno, A., Lissia, M., and Quarati, P.: 1998, Anomalous\par
diffusion modifies solar neutrino fluxes, {\it Physica} {\bf A261}, 359-373.\par
\medskip
\noindent
Lang, K.R.: 1999, {\it Astrophysical Formulae Vol. I (Radiation, Gas\par
Processes and High Energy Astrophysics) and Vol. II (Space, Time, Matter\par
and Cosmology)}, Springer-Verlag, Berlin-Heidelberg.\par
\medskip
\noindent
Mathai, A.M.: 1993, {\it A Handbook of Generalized Special Functions for\par
Statistical and Physical Sciences,} Clarendon Press, Oxford.\par
\medskip
\noindent
Mathai, A.M. and Saxena, R.K.: 1978, {\it The H-Function with Applications\par
in Statistical and Other Disciplines,} Wiley Halsted, New York.\par
\medskip
\noindent
Shaviv, N.J. and Shaviv, G.: 2001, Deriving the electrostatic screening of\par
nuclear reactions from first principles, {\it Nuclear Physics} {\bf A688},\par
285c-288c.\par
\medskip
\noindent
Ueda, M., Sargeant, A.J., Pato, M.P., and Hussein, M.S.: 2000, Effective\par
astrophysical S factor for nonresonant reactions, {\it Physical Review} {\bf C61}\par
045801-1.\par
\medskip
\noindent
Wierling, Th., Millat, and Roepke, G.: 2001, Dynamical screening corrections\par
to the electron capture rate by 7Be, {\it Nuclear Physics} {\bf A688}, 569c-571c.\par
\medskip
\noindent
Wolfram, S.: 1999, {\it The Mathematica Book}, Fourth Edition (Mathematica\par
Version 4), Wolfram Media and Cambridge University Press, Cambridge.\par
\medskip
\noindent
Wolfram Research's Mathematical Functions: 2002, http://functions.wolfram.com/
\end{document}